\documentclass[a4paper,12pt]{article}

 \usepackage{amsfonts,amssymb,mathrsfs}
 \usepackage{times,helvet,courier,type1cm}
 \usepackage{color}
 \usepackage[svgnames]{xcolor}
 \usepackage[leqno]{amsmath}

 \allowdisplaybreaks

 \newtheorem{thm0}{Theorem}[section]
 \newtheorem{def0}{Theorem}[section]

 \newtheorem{exa0}{Theorem}[section]
 
 \newtheorem{thm1}[thm0]{Theorem}
 \newtheorem{def1}[def0]{Definition}
 \newtheorem{lem1}[thm0]{Lemma}
 \newtheorem{cor1}[thm0]{Corollary}
 \newtheorem{pro1}[thm0]{Proposition}
 \newtheorem{rem1}[thm0]{Remark}

 \def\btheorem{\begin{thm1}}\def\etheorem{\end{thm1}}
 \def\bcorollary{\begin{cor1}}\def\ecorollary{\end{cor1}}
 
 \def\bremark{\begin{rem1}}\def\eremark{\end{rem1}}
 \newtheorem{exa1}[exa0]{Example}
 \def\bproof{\begin{proof}}\def\eproof{\qed\end{proof}}

 \def\benumerate{\begin{enumerate}}\def\eenumerate{\end{enumerate}}
 \def\bitemize{\begin{itemize}}\def\eitemize{\end{itemize}}\def\itm{\item}
 \def\bexample{\begin{exa1}}\def\eexample{\end{exa1}}

 \def\beqlb{\begin{eqnarray}}\def\eeqlb{\end{eqnarray}}
 \def\beqnn{\begin{eqnarray*}}\def\eeqnn{\end{eqnarray*}}

 \def\eqref#1{{\rm(\ref{#1})}}

 \def\ar{\!\!\!&}\def\nnm{\nonumber}\def\ccr{\nnm\\}

 \def\<{\langle}\def\>{\rangle}

 \def\mcr{\mathscr}\def\mbb{\mathbb}\def\mbf{\mathbf}
 \def\mrm{\mathrm}

 \def\qed{\hfill$\square$\smallskip}

 \def\d{\mrm{d}}\def\e{\mrm{e}}\def\im{\mrm{i}}

 \def\itDelta{{\it\Delta}}
 
 \def\itOmega{{\it\Omega}}

 \def\qqquad{\qquad\qquad}

\def\rv{}\def\rrv{}

\begin{document}

\noindent{(Version: 2023-02-24)}

\bigskip\bigskip

\centerline{\Large Uniqueness Problem for the Backward Differential}

\smallskip

\centerline{\Large Equation of a Continuous-State Branching Process\footnote{This work was supported by the National Key R\&D Program of China (No.~2020YFA0712900) and the National Natural Science Foundation of China (No.12271029).}}

\bigskip

\centerline{Pei-Sen Li\,\footnote{School of Mathematics and Statistics, Beijing Institute of Technology, Beijing 100872, People's Republic of China (peisenli@bit.edu.cn)} and Zenghu Li\,\footnote{School of Mathematical Sciences, Beijing Normal University, Beijing 100875, People's Republic of China (lizh@bnu.edu.cn)}}

\bigskip\bigskip

\textbf{Abstract:} The distributional properties of a multi-dimensional continuous-state branching process are determined by its cumulant {\rv semigroup, which is defined} by the backward differential equation. We provide a proof of the assertion of Rhyzhov and Skorokhod (Theory Probab. Appl., 1970) on the uniqueness of the solutions to the equation, which is based on a characterization of the process as the pathwise unique solution to a system of stochastic equations.

\bigskip

\smallskip


\bigskip
\section{Introduction}

 \setcounter{equation}{0}

Let $\mbb{R}_+= [0,\infty)$ and let $D= \mbb{R}_+^m$ be its $m$th power for an integer $m\ge 1$. Suppose that $\xi= \{\xi(t): t\ge 0\}$ is a stochastically continuous Markov process in $D$ with homogeneous transition semigroup $\{P(t): t\ge 0\}$. We call $\xi$ a \textit{continuous-state branching process} (CB-process) if its transition semigroup satisfies the \textit{branching property}:
 \beqlb\label{P(t,x+y,.)=}
P(t,x+y,\cdot) = P(t,x,\cdot)*P(t,y,\cdot), \quad x,y\in D,
 \eeqlb
where ``$*$'' denotes the convolution operation. Here $\xi$ may not be conservative, so it is in general a Markov process in the one-point compactification $\bar{D}:= D\cup \{\infty\}$ with $\infty$ being the cemetery.

{\rv The CB-process is a special form of the measure-valued branching processes studied in Dawson \cite{Daw93}, Dynkin \cite{Dyn94}, Le~Gall \cite{LeG99} and Li \cite{Li22}. A more general class of finite-dimensional Markov processes, the \textit{affine processes}, were investigated in Duffie et al.\ \cite{DFS03}. The process has played important roles in the study of different areas including biological populations, statistical physics, mathematical finance and so on. The problem of characterizing the most general branching Markov processes in finite- or infinite-dimensional state spaces has been studied by Duffie et al.\ \cite{DFS03}, Dynkin et al.\ \cite{DKS94}, Rhyzhov and Skorokhod \cite{RhS70}, Silverstein \cite{Sil68} and Watanabe \cite{Wat69}.}

{\rv The basic structures of the multi-dimensional CB-process were investigated by Rhyzhov and Skorokhod \cite{RhS70} and Watanabe \cite{Wat69}.} Let $\mbb{C}_-= \{x+\im y: x\le 0, y\in \mbb{R}\}$, where $\im= \sqrt{-1}$, be the closed left half complex plane and let $\mbb{C}_{--}= \{x+\im y: x< 0, y\in \mbb{R}\}$ be its interior. The $m$th powers of them are denoted by $\mbb{C}_-^m$ and $\mbb{C}_{--}^m$, respectively. From the branching property \eqref{P(t,x+y,.)=}, it follows that $P(t,x,\cdot)$ is an infinitely divisible distribution on $D$. A more careful analysis shows that the Laplace transform of $P(t,x,\cdot)$ is characterized by, for $\lambda\in \mbb{C}_-^m$,
 \beqlb\label{int_De^<>P(t,x,dz)=}
\int_D \e^{\<\lambda,z\>} P(t,x,\d z) = \exp\big\{\<x,K(t,\lambda)\>\big\},
 \eeqlb
where $\<\cdot,\cdot\>$ denotes the inner product and $K(t,\lambda)= (K_1(t,\lambda), \cdots, K_m(t,\lambda))\in \mbb{C}_-^m$. From \eqref{int_De^<>P(t,x,dz)=} we have, for $i= 1,\cdots,m$,
 \beqlb\label{K_i(t,.)=log}
K_i(t,\lambda)= \log\int_D \e^{\<\lambda,z\>} P(t,\delta_i,\d z),
 \eeqlb
where $\delta_i= (1_{\{i=1\}}, \cdots, 1_{\{i=m\}})\in D$. The Chapman--Kolmogorov equation of $\{P(t): t\ge 0\}$ implies
 \beqlb\label{K(s+t,.)= K(s,K(t,.))}
K(s+t,\lambda)= K(s,K(t,\lambda)), \quad s,t\ge 0,\lambda\in \mbb{C}_-^m.
 \eeqlb
We call $\{K(t): t\ge 0\}$ the \textit{cumulant semigroup} of $\xi$ or $\{P(t): t\ge 0\}$. In view of \eqref{int_De^<>P(t,x,dz)=}, we have
 \beqlb\label{P(t,x,{0})=}
P(t,x,D) = \exp\big\{\<x,K(t,0)\>\big\}, \quad t\ge 0,x\in D.
 \eeqlb

From the stochastic continuity of $\{\xi(t): t\ge 0\}$ it follows that $t\mapsto K(t,\lambda)$ is continuous in $t\ge 0$ {\rv for each $\lambda\in \mbb{C}_-^m$ and the restriction of $\lambda\mapsto K(t,\lambda)$ to $\mbb{C}_{--}^m$ is a holomorphic self-map. In fact, the infinite divisibility of the transition probabilities yields a L\'{e}vy--Khintchine type representation of the function $\lambda\mapsto K(t,\lambda)$. For $m=1$ this is known as a \textit{Bernstein function}; see, e.g., Schilling et al.\ \cite{SSV12}. It was shown by Rhyzhov and Skorokhod \cite{RhS70} that for any $\lambda\in \mbb{C}_{--}^m$ the following derivative exists:
 \beqlb\label{H=dK/dt(0,.)}
H(\lambda)= \frac{\d}{\d t} K(t,\lambda)\Big|_{t=0}.
 \eeqlb
Rhyzhov and Skorokhod \cite{RhS70} also derived from \eqref{K(s+t,.)= K(s,K(t,.))} and \eqref{H=dK/dt(0,.)} the differential equation, for $\lambda\in \mbb{C}_{--}^m$,
 \beqlb\label{dK/dt(t)=H(K(t))}
\frac{\d}{\d t} K(t,\lambda)= H(K(t,\lambda)), \quad K(0,\lambda) = \lambda.
 \eeqlb
This is referred to as the \textit{backward differential equation} of the cumulant semigroup since it corresponds to Kolmogorov's backward differential equation of the transition semigroup $\{P(t): t\ge 0\}$. A similar differential equation for the generating functions of integer-valued branching processes was given by Athreya and Ney \cite[p.106]{AtN72}.}

We call $H= (H_1,\cdots,H_m)$ the \textit{branching mechanism} of $\xi$. In complex analysis, one usually refers to \eqref{dK/dt(t)=H(K(t))} as the \textit{Loewner--Kufarev equation}; see, e.g., Bracci et al.\ \cite{BCD20}, Duren \cite{Dur83} and Gumenyuk et al.\ \cite{GHP22a+, GHP22b+}. It is also called the \textit{generalized Riccati equation}; see Duffie et al.\ \cite{DFS03} and the references therein. {\rv For $x= (x_1,\cdots,x_m)\in D$ let $|x|= |x_1| + \cdots + |x_m|$. In Rhyzhov and Skorokhod \cite{RhS70}, the following L\'{e}vy--Khintchine type representation of $H= (H_1,\cdots,H_m)$ was identified: for $i= 1, \cdots, m$,
 \beqlb\label{H_i=<a_i,.>+.}
H_i(\lambda)= \<a_i,\lambda\> + \frac{1}{2}\beta_i\lambda_i^2 + \int_{D^\circ} \Big(\e^{\<\lambda,z\>} - 1 - \frac{\lambda_iz_i}{1+|z|^2}\Big)\pi_i(\d z),
 \eeqlb
where $\beta_i\ge 0$, $a_i= (a_{i1},\cdots,a_{im})\in \mbb{R}^m$ satisfies $a_{ij}\ge 0$ for $j\neq i$ and $\pi_i(\d z)$ is a $\sigma$-finite measure on $D^\circ:= D\setminus\{0\}$ satisfying
 \beqlb\label{int(..1)pi_i(dz)<infty}
\int_{D^\circ} \Big[\Big(\sum_{j=1,j\neq i}^m z_j + z_i^2\Big)\land 1\Big]\pi_i(\d z)< \infty.
 \eeqlb

The property $a_{ij}\ge 0$ for $j\neq i$ follows essentially from \eqref{H=dK/dt(0,.)} and the fact that $\lambda\mapsto K(t,\lambda)$ is the cumulant of an infinitely divisible distribution on $\mbb{R}_+^m$. It was derived rigorously in Duffie et al.\ \cite[Theorem~2.7]{DFS03}. Note that the property was stated incorrectly as $a_{ij}> 0$ for $j\neq i$ in \cite[p.706]{RhS70}. The integrability condition \eqref{int(..1)pi_i(dz)<infty} also corrects the condition given in \cite[p.706]{RhS70}, where it was only claimed that
 \beqlb\label{int_{.}()pi_i(dz)<infty}
\int_{\{|z|<1\}} \Big(\sum_{j=1,j\neq i}^m z_j + z_i^2\Big) \pi_i(\d z)< \infty.
 \eeqlb
It is easy to see that \eqref{int(..1)pi_i(dz)<infty} is satisfied if and only if \eqref{int_{.}()pi_i(dz)<infty} holds and
 \beqlb\label{pi_i(|z|ge1)<infty}
\pi_i(\{z\in D^\circ: |z|\ge 1\})< \infty,
 \eeqlb
which is a standard assumption for the L\'{e}vy measure; see, e.g., Sato \cite[Theorem~8.1]{Sat99}. For instance, if $\lambda= -(1,\cdots,1)$, then the integrand in \eqref{H_i=<a_i,.>+.} becomes
 \beqnn
f(z)= \e^{-|z|} - 1 + \frac{z_i}{1+|z|^2},
 \eeqnn
which satisfies
 \beqnn
-1\le f(z)\le \e^{-|z|} - 1 + \frac{|z|}{1+|z|^2}
 \le
- \Big(\frac{1}{2}-\frac{1}{\e}\Big), \quad |z|\ge 1.
 \eeqnn
Then the integral on the right-hand side of \eqref{H_i=<a_i,.>+.} diverges when \eqref{pi_i(|z|ge1)<infty} is not satisfied. This shows that \eqref{int(..1)pi_i(dz)<infty} is the correct assumption.

For the convenience of our exploration, we need another representation of the branching mechanism $H= (H_1,\cdots,H_m)$. Let $\alpha_i= (\alpha_{i1},\cdots,\alpha_{im})$, where $\alpha_{ij}= a_{ij}$ for $j\neq i$ and
 \beqnn
\alpha_{ii}= a_{ii} + \int_{D^\circ} z_i\Big(1_{\{|z|\le 1\}} - \frac{1}{1+|z|^2}\Big) \pi_i(\d z).
 \eeqnn
Then, for $i= 1, \cdots, m$,
 \beqlb\label{H_i=<alpha_i,.>+.}
H_i(\lambda)= \<\alpha_i,\lambda\> + \frac{1}{2}\beta_i\lambda_i^2 + \int_{D^\circ} \big(\e^{\<\lambda,z\>} - 1 - \lambda_iz_i1_{\{|z|\le 1\}}\big)\pi_i(\d z).
 \eeqlb

It is known that the one-dimensional CB-process is conservative if and only if the following condition is satisfied:
 \beqlb\label{eq1.4}
\int_{0+}\frac{\d\lambda}{0\vee [-H(-\lambda)]}= \infty;
 \eeqlb
see, e.g., \cite{Gre74, IkW70, KaW71}. There seems no necessary and sufficient condition in such a simple form for the conservativeness of multi-dimensional CB-processes.

\bexample\label{ex-alpha-stable} Let us consider the one-dimensional case, i.e., $m=1$. For $\sigma> 0$ and $0<\alpha\le 1$, an \textit{$\alpha$-stable branching mechanism} is defined by $H(\lambda)= -\sigma(-\lambda)^\alpha$. In this case, for $\lambda\in \mbb{C}_{--}$ the unique solution to \eqref{dK/dt(t)=H(K(t))} is given by
 \beqnn\rrv
K(t,\lambda)
 =\bigg\{
\begin{array}{ll}
 -[\sigma(1-\alpha)t+(-\lambda)^{1-\alpha}]^{1/(1-\alpha)} \ar~~ \mbox{for $0<\alpha< 1$,} \cr
 \lambda\e^{\sigma t} \ar~~\mbox{for $\alpha= 1$.}
\end{array}
 \eeqnn
This is a holomorphic function in $\mbb{C}_{--}$ and has a unique continuous extension on $\mbb{C}_-$ with
 \beqnn\rrv
K(t,0)
 =\bigg\{
\begin{array}{ll}
 -[\sigma(1-\alpha)t]^{1/(1-\alpha)} \ar~~ \mbox{for $0<\alpha< 1$,} \cr
 0 \ar~~\mbox{for $\alpha= 1$,}
\end{array}
 \eeqnn
which solves \eqref{dK/dt(t)=H(K(t))} for $\lambda= 0$. The solution to \eqref{dK/dt(t)=H(K(t))} is unique for $\lambda= 0$ if and only if $\alpha= 1$, as $K^\infty(t,0):= 0$ also solves the equation for $\lambda= 0$. \eexample

\bexample\label{ex-(1/2)-stable} Let $m=1$ and consider the $\frac{1}{2}$-stable branching mechanism $H(\lambda)= - 2\sqrt{-\lambda}$. In this case, for $\lambda\in \mbb{C}_{--}$ we have
 \beqnn
K(t,\lambda)= -(t+\sqrt{-\lambda})^2,
 \eeqnn
which maps $\mbb{C}_{--}$ to the domain $\{x+\im y: x< -t\sqrt{t^2+2|y|}\}$. The uniqueness of solution to \eqref{dK/dt(t)=H(K(t))} breaks down for $\lambda= 0$. In fact, for any $r\in [0,\infty]$ the map $t\mapsto K^r(t,0):= -(t-r)^21_{\{t> r\}}$ is a solution to \eqref{dK/dt(t)=H(K(t))} for $\lambda= 0$. \eexample

\bremark\label{th-int-eqK(t,.)} Let $H$ be a branching mechanism defined by \eqref{H_i=<a_i,.>+.} or \eqref{H_i=<alpha_i,.>+.}. From the result of Duffie et al.\ \cite[Theorem~2.7 and Proposition~6.1]{DFS03} it follows that for any $\lambda\in \mbb{C}_{--}^m$ there is a unique solution to \eqref{dK/dt(t)=H(K(t))} and one can define the transition semigroup of a CB-process by \eqref{int_De^<>P(t,x,dz)=}; see also Watanabe~\cite[Theorems~2 and~2$'$]{Wat69}. By \eqref{K_i(t,.)=log} we can extend $\lambda\mapsto K(t,\lambda)$ to a continuous self-map of $\mbb{C}_-^m$. For $\lambda\in \mbb{C}_{--}^m$ we can rewrite \eqref{dK/dt(t)=H(K(t))} into its integral form
 \beqlb\label{K(t)=.+int_0^tH()ds}
K(t,\lambda)= \lambda + \int_0^tH(K(s,\lambda))\d s, \quad t\ge 0.
 \eeqlb
In view of \eqref{H_i=<a_i,.>+.} or \eqref{H_i=<alpha_i,.>+.}, it is clear that the function $\lambda\mapsto H(\lambda)$ admits a unique continuous extension on $\mbb{C}_-^m$ with $H(0)= 0$. By approximation from the interior and application of dominated convergence we see that \eqref{K(t)=.+int_0^tH()ds} remains true for $\lambda\in \partial \mbb{C}_-^m:= \mbb{C}_-^m\setminus \mbb{C}_{--}^m$. By differentiating both sides of the integral equation we see that \eqref{dK/dt(t)=H(K(t))} also holds for all $\lambda\in \mbb{C}_-^m$. That proves the existence of the solution to \eqref{dK/dt(t)=H(K(t))} for all $\lambda\in \mbb{C}_-^m$; see also Duffie et al.\ \cite[Proposition~6.4 and its proof]{DFS03}. \eremark

Rhyzhov and Skorokhod \cite[p.706]{RhS70} wrote: It can be shown that if \eqref{dK/dt(t)=H(K(t))} has a unique solution belonging to $\mbb{C}_-^m$ for all $\lambda\in \mbb{C}_-^m$, then $\exp\{\langle x,K(t,\lambda)\rangle\}$ will be a transform of type \eqref{P(t,x+y,.)=} of the transition probability of a branching process. The following assertion is very plausible: \textit{in order that the equation \eqref{dK/dt(t)=H(K(t))} have a unique solution for all $\lambda\in \mbb{C}_-^m$, it is necessary and sufficient that it have a unique solution for $\lambda= 0$.} In view of Remark~\ref{th-int-eqK(t,.)}, we can restate their assertion equivalently as: in order that the equation \eqref{dK/dt(t)=H(K(t))} have a unique solution for all $\lambda\in \mbb{C}_-^m$, it is necessary and sufficient that $t\mapsto K(t,0)\equiv 0$ is the unique solution to \eqref{dK/dt(t)=H(K(t))} for $\lambda= 0$.

}

In this paper, we provide a proof of the assertion of Rhyzhov and Skorokhod \cite[p.706]{RhS70}. The proof is based on a characterization of the $m$-dimensional CB-process as the pathwise unique solution to a system of stochastic equations.

\section{The system of stochastic equations}

 \setcounter{equation}{0}

Recall that the vector-valued function {\rv $H= (H_1,\cdots,H_m)$} is defined by \eqref{H_i=<alpha_i,.>+.} with parameters $(a_i,\beta_i,\pi_i)$, $i=1,\cdots,m$. Suppose that $(\itOmega,\mcr{F},\mcr{F}_t,\mbf{P})$ is a filtered probability space satisfying the usual hypotheses. For each $k= 1, \cdots, m$ let $\{B_k(t)\}$ be an $(\mcr{F}_t)$-Brownian motion and let $\{N_k(\d s,\d z,\d u)\}$ be a time-space $(\mcr{F}_t)$-Poisson random measure on $(0,\infty)\times D^\circ\times (0,\infty)$ with intensity $\d s\pi_k(\d z)\d u$. Let $\tilde{N}_k(\d s,\d z,\d u):= N_k(\d s,\d z,\d u) - \d s\pi_k(\d z)\d u$ denote the compensated measure. Let $D_1= \{x\in D^\circ: |x|\le 1\}$ and $D_1^c= \{x\in D^\circ: |x|> 1\}$. We consider the system of stochastic integral equations, for $k= 1, \cdots, m$,
 \beqlb\label{xi_k(t)=xi_k(0)+.}
\xi_k(t)\ar=\ar \xi_k(0) + \int_0^t\sqrt{\beta_k \xi_k(s)} \d B_k(s) + \sum_{i=1}^m\int_0^t \alpha_{ik}\xi_i(s) \d s \cr
 \ar\ar
+ \int_0^t\int_{D_1}\int_0^{\xi_k(s-)} z_k\tilde{N}_k(\d s,\d z,\d u) + \int_0^t\int_{D_1^c}\int_0^{\xi_k(s-)} z_kN_k(\d s,\d z,\d u) \cr
 \ar\ar
+ \sum_{i=1,i\neq k}^m\int_0^t\int_{D^\circ}\int_0^{\xi_i(s-)} z_kN_i(\d s,\d z,\d u).
 \eeqlb
We say an $(\mcr{F}_t)$-adapted c\`{a}dl\`{a}g process $\{\xi(t): t\ge 0\}= \{(\xi_1(t),\cdots,\xi_m(t)): t\ge 0\}$ {\rv in $D$ with possibly finite lifetime is a \textit{solution} to \eqref{xi_k(t)=xi_k(0)+.} if the equation holds almost surely when $t$ is replaced by $t\land \zeta_n$ for every $t\ge 0$ and $n\ge 1$, where $\zeta_n= \inf\{t\ge 0: |\xi(t)|\ge n\}$ with $\inf\emptyset = \infty$ by convention. We make the convention that $\xi(t)= \infty$ for $t\ge \zeta:= \lim_{n\to \infty} \zeta_n$. By saying $\{\xi(t)\}$ is \textit{conservative} we mean $\zeta= \infty$ almost surely.

To prove the existence and uniqueness of the solution to \eqref{xi_k(t)=xi_k(0)+.}, we introduce an approximating sequence of the stochastic equation. For each $n\ge 1$ let us consider the equations, for $k= 1, \cdots, m$,
 \beqlb\label{x_{nk}(t)=xi_k(0)+.}
x_{nk}(t)\ar=\ar \xi_k(0) + \int_0^t\sqrt{\beta_k x_{nk}(s)} \d B_k(s) + \sum_{i=1}^m\int_0^t \alpha_{ik}x_{ni}(s) \d s \cr
 \ar\ar
+ \int_0^t\int_{D_1}\int_0^{x_{nk}(s-)} z_k\tilde{N}_k(\d s,\d z,\d u) + \int_0^t\int_{D_1^c}\int_0^{x_{nk}(s-)} (z_k\land n)N_k(\d s,\d z,\d u) \quad \cr
 \ar\ar
+ \sum_{i=1,i\neq k}^m\int_0^t\int_{D^\circ}\int_0^{x_{ni}(s-)} (z_k\land n)N_i(\d s,\d z,\d u),
 \eeqlb
The truncation $z_k\land n:= \min\{z_k, n\}$ in \eqref{x_{nk}(t)=xi_k(0)+.} changes the jump sizes larger than $n$ into jump size $n$. With this modification of the equation, we can apply Barczy et al.\ \cite[Theorem~4.6]{BLP15} to see that there is a pathwise unique solution $\{x_n(t):t\ge 0\}= \{(x_{n1}(t),\cdots,x_{nm}(t)): t\ge 0\}$ to \eqref{x_{nk}(t)=xi_k(0)+.}.}

\btheorem\label{th-path-uniq} For any $\mcr{F}_0$-measurable random variable $\xi(0)\in D$, there is a pathwise unique solution $\{\xi(t): t\ge 0\}$ to \eqref{xi_k(t)=xi_k(0)+.}. \etheorem

\bproof {\rv For $n\ge 1$ define the stopping time $\zeta_n= \inf\{t\ge 0: |x_n(t)|\ge n\}$. Then $\rrv |x_n(t)|< n$ for $0\le t< \zeta_n$. Since the processes $\{x_{nk}(t): t\ge 0\}$, $k= 1,\cdots,m$ are nonnegative and can only have positive jumps, they do not have jumps with jump sizes larger than $n$ in the time interval $(0,\zeta_n)$, so the truncation $(z_k\land n)$ in \eqref{x_{nk}(t)=xi_k(0)+.} makes no difference in this time interval. It follows that $x_{n+1}(t)= x_n(t)$ for $0\le t< \zeta_n$ and $\zeta_n$ is nondecreasing in $n\ge 1$. Then we can define a process $\{\xi(t): t\ge 0\}= \{(\xi_1(t),\cdots,\xi_m(t)): t\ge 0\}$ in $\bar{D}$ such that $\xi(t)= x_n(t)$ for $0\le t< \zeta_n$ and $\xi(t)= \infty$ for $t\ge \zeta:= \lim_{n\to \infty} \zeta_n$. It is easy to see that $\{\xi(t): t\ge 0\}$ is a solution to \eqref{xi_k(t)=xi_k(0)+.} and $\zeta_n= \inf\{t\ge 0: |\xi(t)|\ge n\}$. This proves the existence of solution to \eqref{xi_k(t)=xi_k(0)+.}. Now suppose that $\{y(t): t\ge 0\}= \{(y_1(t),\cdots,y_m(t)): t\ge 0\}$ is also a solution to \eqref{xi_k(t)=xi_k(0)+.}. Let $\tau_n= \inf\{t\ge 0: |y(t)|\ge n\}$ and $\eta_n= \tau_n\land \zeta_n$ for $n\ge 1$. Then $|y(t)|< n$ and $|\xi(t)|< n$ for $0\le t< \eta_n$, and so the processes $\{y(t): t\ge 0\}$ and $\{\xi(t): t\ge 0\}$ coincide in the time interval $(0,\eta_n)$. In view of \eqref{xi_k(t)=xi_k(0)+.} and \eqref{x_{nk}(t)=xi_k(0)+.} we have $\xi(\eta_n)= y(\eta_n)$, which implies $\eta_n= \zeta_n= \tau_n$. It follows that $\zeta= \lim_{n\to \infty} \zeta_n= \lim_{n\to \infty} \eta_n$, and so the processes $\{y(t): t\ge 0\}$ and $\{\xi(t): t\ge 0\}$ are indistinguishable. Then the pathwise uniqueness of solutions holds for \eqref{xi_k(t)=xi_k(0)+.}.} \eproof

The system of stochastic equations \eqref{xi_k(t)=xi_k(0)+.} gives a construction of the sample path of the $m$-dimensional CB-process. The existence and pathwise uniqueness of solutions to \eqref{xi_k(t)=xi_k(0)+.} were established in Dawson and Li \cite{DaL06} for $m=1$ under a stronger moment condition. A flow of discontinuous CB-processes was constructed in Bertoin and Le~Gall \cite{BeL06} by weak solutions to a special form of the one-dimensional stochastic equation; see also \cite{BeL00}. Their construction was extended to general flows in Dawson and Li \cite{DaL12} using strong solutions. See, e.g., Bernis et al.\ \cite{BBSS21}, Fang and Li \cite{FaL22}, Li \cite{Li20, Li22} and Pardoux \cite{Par16} for a number of applications of those stochastic equations. {\rv Theorem~\ref{th-path-uniq} weakens the moment conditions in Barczy et al.\ \cite{BLP15} and Ma \cite{Ma13}, where multi-dimensional continuous state branching processes with immigration were constructed in terms of stochastic equations.}

\section{A time-space martingale problem}

 \setcounter{equation}{0}

Let $C_0^2(D)$ be the set of twice continuously differentiable functions $f= f(x)$ on $D$ that together with their derivatives up to the second order are rapidly decreasing as $|x|\to \infty$. We extend the functions $f= f(x)$ in $C_0^2(D)$ and those derivatives trivially to $\infty$ so that they become bounded continuous functions on $\bar{D}$. {\rv For $x\in D$ and $f\in C_0^2(D)$ set
 \beqlb\label{Af(x)=...}
Af(x)\ar=\ar \sum_{k=1}^m x_k\<\alpha_k,\nabla f(x)\> + \frac{1}{2}\sum_{k=1}^m\beta_k x_k f_{kk}''(x) \cr
 \ar\ar
+ \sum_{k=1}^m x_k\int_{D^\circ} \big[\itDelta_zf(x)-z_kf'_k(x)1_{D_1}(z)\big] \pi_k(\d z),
 \eeqlb
where $\itDelta_zf(x)= f(x+z)-f(x)$. Then we can understand $Af$ as a bounded continuous function on $\bar{D}$ with $Af(\infty)= 0$.

For a fixed constant $u>0$ let $C^{1,2}([0,u]\times D)$ be the set of functions $f= f(t,x)$ on $[0,u]\times D$ that are continuously differentiable in $t\in [0,u]$ and twice continuously differentiable in $x= (x_1,\cdots,x_m)\in D$. For $f\in C^{1,2}([0,u]\times D)$ and $i,j=1,\cdots,m$ write
 \beqnn
f_0'(t,x)= \frac{\partial}{\partial t} f(t,x),~ f_i'(t,x)= \frac{\partial}{\partial x_i} f(t,x),~ f_{ij}''(t,x)= \frac{\partial}{\partial x_ix_j} f(t,x).
 \eeqnn
Let $C_0^{1,2}([0,u]\times D)$ be the set of functions $f= f(t,x)$ in $C^{1,2}([0,u]\times D)$ that together with the above derivatives are rapidly decreasing as $|x|\to \infty$. We extend all the functions in $C_0^{1,2}([0,u]\times D)$ and those derivatives trivially to $[0,u]\times \{\infty\}$ so that they become bounded continuous functions on $[0,u]\times \bar{D}$.} The theorem below gives the characterization of the solution $\{\xi(t): t\ge 0\}$ to \eqref{xi_k(t)=xi_k(0)+.} by a time-space martingale problem, which identifies the operator $A$ defined by \eqref{Af(x)=...} as (a restriction of) its generator.

\btheorem\label{time-space-mart} Let $\{\xi(t): t\ge 0\}$ be the pathwise unique solution to \eqref{xi_k(t)=xi_k(0)+.}. {\rv Then {\rrv for any} $0\le t\le u$ and $f\in C_0^{1,2}([0,u]\times D)$ we have}
 \beqlb\label{f(t,xi(t))=f(0,xi(0))+.}
f(t,\xi(t))= f(0,\xi(0)) + \int_0^t \big[f_0'(s,\xi(s))+Af(s,\xi(s))\big] \d s + M(t),
 \eeqlb
{\rv where $A$ acts on the function $x\mapsto f(s,x)$ and $\{M(t): 0\le t\le u\}$ is a bounded martingale defined by}
 \beqlb\label{M(t)=...}\rv
M(t)\ar\rv=\ar\rv \sum_{k=1}^m \int_0^t\sqrt{\beta_k \xi_k(s)}f_k'(s,\xi(s)) \d B_k(s) \cr
 \ar\ar\rv
+ \sum_{k=1}^m\int_0^t\int_{D^\circ}\int_0^{\xi_k(s-)} \itDelta_zf(s,\xi(s-)) \tilde{N}_k(\d s,\d z,\d u).
 \eeqlb
\etheorem

\bproof For $n\ge 1$ let $\{x_n(t):t\ge 0\}= \{(x_{n1}(t),\cdots,x_{nm}(t)): t\ge 0\}$ be the pathwise unique solution to \eqref{x_{nk}(t)=xi_k(0)+.}. By It\^o's formula, {\rv for $0\le t\le u$ we have
 \beqnn
\ar\ar f(t,x_n(t)) = f(0,\xi(0)) + \int_0^t f_0'(s,x_n(s-)) \d s + \sum_{k=1}^m\sum_{i=1}^m\int_0^t \alpha_{ik}x_{ni}(s-)f_k'(s,x_n(s-)) \d s \cr
 \ar\ar\qqquad
+ \sum_{k=1}^m\int_0^t\sqrt{\beta_k x_{nk}(s-)}f_k'(s,x_n(s-)) \d B_k(s) + \frac{1}{2}\sum_{k=1}^m\int_0^t \beta_k x_{nk}(s-)f_{kk}''(s,x_n(s-)) \d s \cr
 \ar\ar\qqquad
+ \sum_{k=1}^m\int_0^t\int_{D_1}\int_0^{x_{nk}(s-)} z_kf_k'(s,x_n(s-)) \tilde{N}_k(\d s,\d z,\d u) \cr
 \ar\ar\qqquad
+ \sum_{k=1}^m\int_0^t\int_{D_1^c}\int_0^{x_{nk}(s-)} (z_k\land n) f_k'(s,x_n(s-)) N_k(\d s,\d z,\d u) \cr
 \ar\ar\qqquad
+ \sum_{k=1}^m\sum_{i=1,i\neq k}^m\int_0^t\int_{D^\circ}\int_0^{x_{ni}(s-)} (z_k\land n) f_k'(s,x_n(s-)) N_i(\d s,\d z,\d u) \cr
 \ar\ar\qqquad
+ \sum_{k=1}^m\int_0^t\int_{D^\circ}\int_0^{x_{nk}(s-)} \big[f(s,x_n(s-)+(z\land n))-f(s,x_n(s-)) \cr
 \ar\ar\qqquad\qqquad\qqquad\qqquad\qquad
-\, \<z\land n,\nabla f(s,x_n(s-))\>\big] N_k(\d s,\d z,\d u) \ccr
 \ar\ar\qquad\quad
= f(0,\xi(0)) + \int_0^t f_0'(s,x_n(s-)) \d s + \sum_{i=1}^m\int_0^t x_{ni}(s-)\<\alpha_i,\nabla f(s,x_n(s-))\> \d s \cr
 \ar\ar\qqquad
+ \sum_{k=1}^m\int_0^t\sqrt{\beta_k x_{nk}(s-)}f_k'(s,x_n(s-)) \d B_k(s) + \frac{1}{2}\sum_{k=1}^m \int_0^t \beta_k x_{nk}(s-)f_{kk}''(s,x_n(s-)) \d s \cr
 \ar\ar\qqquad
+ \sum_{k=1}^m\int_0^t\int_{D_1}\int_0^{x_{nk}(s-)} z_kf_k'(s,x_n(s-)) \tilde{N}_k(\d s,\d z,\d u) \cr
 \ar\ar\qqquad
+ \sum_{k=1}^m\int_0^t\int_{D^\circ}\int_0^{x_{nk}(s-)} \big[f(s,x_n(s-)+(z\land n)) - f(s,x_n(s-)) \cr
 \ar\ar\qqquad\qqquad\qqquad\qqquad\qquad
-\,z_kf_k'(s,x_n(s-))1_{D_1}(z)\big] N_k(\d s,\d z,\d u) \ccr
 \ar\ar\qquad\quad
= f(0,\xi(0)) + \int_0^t f_0'(s,x_n(s-)) \d s + \sum_{i=1}^m\int_0^t x_{ni}(s-)\<\alpha_i,\nabla f(s,x_n(s-))\> \d s \cr
 \ar\ar\qqquad
+ \sum_{k=1}^m\int_0^t\sqrt{\beta_k x_{nk}(s-)}f_k'(s,x_n(s-)) \d B_k(s) + \frac{1}{2}\sum_{k=1}^m \int_0^t \beta_k x_{nk}(s-)f_{kk}''(s,x_n(s-)) \d s \cr
 \ar\ar\qqquad
+ \sum_{k=1}^m\int_0^t\int_{D^\circ}\int_0^{x_{nk}(s-)} \itDelta_{z\land n}f(s,x_n(s-)) \tilde{N}_k(\d s,\d z,\d u) \cr
 \ar\ar\qqquad
+ \sum_{k=1}^m\int_0^t x_{nk}(s-)\d s\int_{D^\circ} \big[\itDelta_{z\land n}f(s,x_n(s-)) - z_kf_k'(s,x_n(s-))1_{D_1}(z)\big] \pi_k(\d z),
 \eeqnn
where $z\land n= (z_1\land n, \cdots, z_m\land n)$. Then
 \beqnn
f(t,x_n(t))= f(0,\xi(0)) + \int_0^t \big[f_0'(s,x_n(s-))+A_nf(s,x_n(s-))\big] \d s + M_n(t),
 \eeqnn
where
 \beqnn
M_n(t)\ar=\ar \sum_{k=1}^m \int_0^t\sqrt{\beta_k x_{nk}(s-)}f_k'(s,x_n(s-)) \d B_k(s) \cr
 \ar\ar
+ \sum_{k=1}^m \int_0^t\int_{D^\circ}\int_0^{x_{nk}(s-)} \itDelta_{z\land n}f(s,x_n(s-)) \tilde{N}_k(\d s,\d z,\d u)
 \eeqnn
and
 \beqnn
A_nf(s,x)\ar=\ar \sum_{k=1}^mx_k\<\alpha_k,\nabla f(s,x)\> + \frac{1}{2}\sum_{k=1}^m \beta_k x_kf_{kk}''(s,x) \cr
 \ar\ar
+ \sum_{k=1}^m x_k\int_{D^\circ} \big[\itDelta_{z\land n}f(s,x) - z_kf_k'(s,x)1_{D_1}(z)\big] \pi_k(\d z).
 \eeqnn
Let $\zeta_n= \inf\{t\ge 0: |x_n(t)|\ge n\}$ for $n\ge 1$. By the proof of Theorem~\ref{th-path-uniq}, the sequence $\{\zeta_n\}$ is nondecreasing. Moreover, we have $\xi(t)= x_n(t)$ for $0\le t< \zeta_n$ and $\xi(t)= \infty$ for $t\ge \zeta:= \lim_{n\to \infty} \zeta_n$. It follows that
 \beqlb\label{f(t,x_n(t))=f(0,xi(0))+.}
f(t,x_n(t\land \zeta_n))= f(0,\xi(0)) + \int_0^{t\land \zeta_n} \big[f_0'(s,\xi(s-))+A_nf(s,\xi(s-))\big] \d s + M_n(t\land \zeta_n),
 \eeqlb
where
 \beqlb\label{M_n(t)=...}
M_n(t\land \zeta_n)\ar=\ar \sum_{k=1}^m \int_0^{t\land \zeta_n}\sqrt{\beta_k \xi_k(s-)}f_k'(s,\xi(s-)) \d B_k(s) \cr
 \ar\ar
+ \sum_{k=1}^m \int_0^{t\land \zeta_n}\int_{D^\circ}\int_0^{\xi_k(s-)} \itDelta_{z\land n}f(s,\xi(s-)) \tilde{N}_k(\d s,\d z,\d u).
 \eeqlb
The definition of $\zeta_n$ implies that $\lim_{n\to \infty} x_n(\zeta_n)= \infty$. By considering the cases $\zeta\le t$ and $\zeta> t$ separately, we see that $\lim_{n\to \infty} x_n(t\land \zeta_n)= \xi(t\land \zeta)$. Then letting $n\to \infty$ in \eqref{f(t,x_n(t))=f(0,xi(0))+.} and \eqref{M_n(t)=...} we obtain
 \beqnn
f(t,\xi(t\land \zeta))= f(0,\xi(0)) + \int_0^{t\land \zeta} \big[f_0'(s,\xi(s-))+Af(s,\xi(s-))\big] \d s + M_\zeta(t),
 \eeqnn
where
 \beqnn
M_\zeta(t)\ar=\ar \sum_{k=1}^m \int_0^{t\land \zeta}\sqrt{\beta_k \xi_k(s-)}f_k'(s,\xi(s-)) \d B_k(s) \cr
 \ar\ar
+ \sum_{k=1}^m \int_0^{t\land \zeta}\int_{D^\circ}\int_0^{\xi_k(s-)} \itDelta_zf(s,\xi(s-)) \tilde{N}_k(\d s,\d z,\d u).
 \eeqnn
Recall that we understand $f$, $f_k'$ and $Af$ as bounded continuous functions on $[0,u]\times \bar{D}$ with $f(s,\infty)= f_k'(s,\infty)= Af(s,\infty)= 0$ for all $s\in [0,u]$. Then we have \eqref{f(t,xi(t))=f(0,xi(0))+.} and \eqref{M(t)=...}. From (3.2) it is clear that
 \beqnn
|M(t)|\le 2\|f\| + (\|f_0'\|+\|Af\|)u, \quad 0\le t\le u.
 \eeqnn
where $\|\cdot\|$ denotes the supremum norm of functions on $[0,u]\times \bar{D}$. Then $\{M(t): 0\le t\le u\}$ is a bounded martingale by \eqref{M(t)=...}.} \eproof

{\rv

\bcorollary\label{th-sol-cadlag} The pathwise unique solution $\{\xi(t): t\ge 0\}$ to \eqref{xi_k(t)=xi_k(0)+.} is a c\`{a}dl\`{a}g process in $\bar{D}$. \ecorollary

\bproof By applying Theorem~\ref{time-space-mart} to the function $f(t,x)\equiv \e^{-x}$ and $u= 1,2,\cdots$ we see that $\{\e^{-\xi(t)}: t\ge 0\}$ is a c\`{a}dl\`{a}g process in $[0,1]$. Then $\{\xi(t): t\ge 0\}$ is a c\`{a}dl\`{a}g process in $\bar{D}$. \eproof

}

\bcorollary\label{th-m-dim-CBp} The {\rv pathwise unique solution} $\{\xi(t): t\ge 0\}$ to \eqref{xi_k(t)=xi_k(0)+.} is an $m$-dimensional CB-process with transition semigroup $\{P(t): t\ge 0\}$ defined by \eqref{int_De^<>P(t,x,dz)=} and \eqref{dK/dt(t)=H(K(t))}. \ecorollary

\bproof {\rv For $\lambda\in \mbb{C}_{--}^m$ let $t\mapsto K(t,\lambda)\in \mbb{C}_{--}^m$ be the unique solution to \eqref{dK/dt(t)=H(K(t))}. Fix $u> 0$ and define
 \beqnn
f(t,x)= \e^{\<K(u-t,\lambda),x\>}, \quad 0\le t\le u, x\in D.
 \eeqnn
Then $f\in C_0^{1,2}([0,u]\times D)$. It is easy to show that for $k=1,\cdots,m$,
 \beqnn
f_k'(t,x)= K_k(u-t,\lambda)\e^{\<K(u-t,\lambda),x\>},
 ~
f_{kk}''(t,x)= K_k(u-t,\lambda)^2\e^{\<K(u-t,\lambda),x\>}
 \eeqnn
and
 \beqnn
\itDelta_zf(x)= (\e^{\<K(u-t,\lambda),z\>}-1)\e^{\<K(u-t,\lambda),x\>}.
 \eeqnn
By \eqref{dK/dt(t)=H(K(t))} and \eqref{Af(x)=...} we have
 \beqnn
f_0'(t,x)= -\e^{\<K(u-t,\lambda),x\>}\<H(K(u-t,\lambda)),x\>= -Af(t,x).
 \eeqnn
Note that the equalities above extend trivially to $x\in \bar{D}$. An application of Theorem~\ref{time-space-mart} shows that
 \beqnn
\e^{\<K(u-t,\lambda),\xi(t)\>}
 =
\e^{\<K(u,\lambda),\xi(0)\>} + M(t), \quad 0\le t\le u,
 \eeqnn
where $\{M(t): 0\le t\le u\}$ is a bounded martingale. It follows that
 \beqnn
\mbf{E}\big[\e^{\<\lambda,\xi(u)\>}|\mcr{F}_t\big]
 \ar=\ar
\mbf{E}\big[\e^{\<K(u,\lambda),\xi(0)\>} + M(u)|\mcr{F}_t\big] \ccr
 \ar=\ar
\e^{\<K(u,\lambda),\xi(0)\>} + M(t)
 =
\e^{\<K(u-t,\lambda),\xi(t)\>}.
 \eeqnn
Then} $\{\xi(t): t\ge 0\}$ is an $m$-dimensional CB-process with transition semigroup $\{P(t): t\ge 0\}$ defined by \eqref{int_De^<>P(t,x,dz)=} and \eqref{dK/dt(t)=H(K(t))}. \eproof

{\rv

\bremark In general, the process $\{\xi(t): t\ge 0\}$ lives in the extended state space $\bar{D}$. For this reason, we need the functions $f$ and $Af$ to have continuous extensions on $[0,u]\times \bar{D}$ in Theorem~\ref{time-space-mart} and in the proof of Corollary~\ref{th-m-dim-CBp}. This is the reason that we assume $\lambda\in \mbb{C}_{--}^m$ in the proof. \eremark

}

\section{The backward differential equation}

 \setcounter{equation}{0}

{\rv By Remark~\ref{th-int-eqK(t,.)}, there exists the solution $t\mapsto K(t,\lambda)$ to \eqref{dK/dt(t)=H(K(t))} for every $\lambda\in \mbb{C}_-^m$. The next theorem confirms the assertion of Rhyzhov and Skorokhod \cite[p.706]{RhS70}:

\btheorem\label{th-RhSk-assert} The following statements are equivalent:
 \bitemize

\itm[(i)] For every $\lambda\in \mbb{C}_-^m$ there is a unique solution $t\mapsto K(t,\lambda)\in \mbb{C}_-^m$ to \eqref{dK/dt(t)=H(K(t))}.

\itm[(ii)] For $\lambda= 0$ the unique solution to \eqref{dK/dt(t)=H(K(t))} is given by $t\mapsto K(t,0)\equiv 0$.

\itm[(iii)] The CB-process with transition semigroup given by \eqref{int_De^<>P(t,x,dz)=} and \eqref{dK/dt(t)=H(K(t))} is conservative.

 \eitemize
\etheorem

\bproof (i)$\Rightarrow$(ii) This is trivial.

(ii)$\Rightarrow$(iii) Suppose that $t\mapsto K(t,0)\equiv 0$ is the unique solution to \eqref{dK/dt(t)=H(K(t))} for $\lambda= 0$. Then \eqref{P(t,x,{0})=} implies $P(t,x,D)= 1$ for all $t\ge 0$ and $x\in D$, which means the CB-process is conservative.

(iii)$\Rightarrow$(i) Suppose that the CB-process with transition semigroup given by \eqref{int_De^<>P(t,x,dz)=} and \eqref{dK/dt(t)=H(K(t))} is conservative. Let $\{\xi(t): t\ge 0\}$ be the realization of the process defined by \eqref{xi_k(t)=xi_k(0)+.} with $\xi(0)= x\in D$. Then $\{\xi(t): t\ge 0\}$ is a c\`{a}dl\`{a}g process living in $D$. Suppose that $t\mapsto J(t,\lambda)$ is also a solution to \eqref{dK/dt(t)=H(K(t))} for $\lambda\in \mbb{C}_-^m$. Fix $u> 0$ and define
 \beqlb\label{g(t,x)=e^<J(u-t),x>}
g(t,x)= \e^{\<J(u-t,\lambda),x\>}, \quad 0\le t\le u, x\in D.
 \eeqlb
It is easy to see that $g\in C^{1,2}([0,u]\times D)$. By \eqref{xi_k(t)=xi_k(0)+.} and It\^o's formula, as in the proof of Theorem~\ref{time-space-mart} we get
 \beqlb\label{g(t,xi(t))=g(0,xi(0))+.}
\ar\ar g(t,\xi(t))= g(0,\xi(0)) + \int_0^t g_0'(s,\xi(s-)) \d s + \sum_{i=1}^m\int_0^t \xi_i(s-)\<\alpha_i,\nabla g(s,\xi(s-))\> \d s \cr
 \ar\ar\qqquad
+ \frac{1}{2}\sum_{k=1}^m \int_0^t \beta_k \xi_k(s-)g_{kk}''(s,\xi(s-)) \d s + Z(t) \cr
 \ar\ar\qqquad
+ \sum_{k=1}^m\int_0^t \xi_k(s-)\d s\int_{D^\circ} \big[\itDelta_zg(s,\xi(s-)) - z_kg_k'(s,\xi(s-))1_{D_1}(z)\big] \pi_k(\d z),
 \eeqlb
where
 \beqnn
Z(t)\ar=\ar \sum_{k=1}^m \int_0^t\sqrt{\beta_k \xi_k(s-)}g_k'(s,\xi(s-)) \d B_k(s) \cr
 \ar\ar
+ \sum_{k=1}^m \int_0^t\int_{D^\circ}\int_0^{\xi_k(s-)} \itDelta_zg(s,\xi(s-)) \tilde{N}_k(\d s,\d z,\d u).
 \eeqnn
Since $t\mapsto J(t,\lambda)$ solves \eqref{dK/dt(t)=H(K(t))}, one can use \eqref{g(t,x)=e^<J(u-t),x>} to see
 \beqnn
g_0'(s,x)\ar=\ar -\sum_{i=1}^m x_i\<\alpha_i,\nabla g(s,x)\> - \frac{1}{2}\sum_{k=1}^m \beta_k x_kg_{kk}''(s,x) \cr
 \ar\ar
- \sum_{k=1}^m x_k\int_{D^\circ} \big[\itDelta_zg(s,x) - z_kg_k'(s,x)1_{D_1}(z)\big] \pi_k(\d z),
 \eeqnn
From \eqref{g(t,xi(t))=g(0,xi(0))+.} it follows that
 \beqlb\label{e^<J(u-t,.),x(t)>=.}
\e^{\<J(u-t,\lambda),\xi(t)\>}= \e^{\<J(u,\lambda),x\>} + Z(t).
 \eeqlb
Then $\{Z(t): 0\le t\le u\}$ is a bounded martingale. (Here we cannot apply Theorem~\ref{time-space-mart} directly since $g\notin C_0^{1,2}([0,u]\times D)$.) By taking the expectations in both sides of \eqref{e^<J(u-t,.),x(t)>=.} for $t=u$ we obtain $\mbf{E}[\e^{\<\lambda,\xi(u)\>}]= \e^{\<J(u,\lambda),x\>}$. Then $J(u,\lambda)= K(u,\lambda)$ by \eqref{int_De^<>P(t,x,dz)=} and Corollary~\ref{th-m-dim-CBp}. This shows the implication ``(iii)$\Rightarrow$(i)''. \eproof

It would be interesting to have an analytic proof of the assertion of Rhyzhov and Skorokhod \cite[p.706]{RhS70}.

\textbf{Acknowledgements.} We would like to thank the two anonymous referees for their valuable comments and suggestions, which have led to a number of improvements in the presentation of the work.

}

\end{document}